\documentstyle[12pt]{article}
\date{}
\newcommand{\vs}{\vskip2mm}

\newcommand{\VS}{\vskip6mm}
\newcounter{tma}
\renewcommand{\thetma}{\arabic{tma}.}
                      {\par\vs}
\newcounter{rem}
\renewcommand{\therem}{\arabic{rem}.}
                   {\par\vs}
\newcounter{gp}
\renewcommand{\thegp}{\arabic{gp}.}
                      {\par}
                         {\par\vs}
\begin{document}
%%%%author
\centerline{\large\bf L.~A.~Sakhnovich}

%%%%title
\begin{center}
      {\Large\bf

     Extremal Trigonometrical and Power Polynomials of Several
     Variables}
\end{center}
\VS
%%%%%%%%%%

%%%%%%%%%%%%%%%%%%%%%%%%%%%%%%%%%%%%%%%%%%%%%%%%%%%%%%%%%%%%%%%%%%%
\setcounter{section}{-1}
\section{\large Introduction.} \,
We consider the set ${\sigma}_{P}$ of the power non-negative
polynomials of several variables.By $Q_{P}$ we denote the class of
the polynomials from  ${\sigma}_{1}$ which can be represented as a
sum of squares.It is shown in the classic work by D.Hilbert[3]
that $Q_{P}$ does not coincide with ${\sigma}_{P}$.Step by step a
number of polynomials belonging to ${\sigma}_{P}$ but not
belonging to $Q_{P}$ was constructed(see[4]-[6]).It is interesting
to note that many of these polynomials turn to be extremal in the
class  ${\sigma}_{P}$ [2].\\
In our paper we have made an attempt to work out a general
approach to the investigation of the extremal elements of the
convex sets $Q_{P}$ and ${\sigma}_{P}$.It seems to us that we have
achieved a considerable progress in the case of $Q_{P}$.In the
case of ${\sigma}_{P}$ we have made only the first steps.We also
consider the class ${\sigma}_{R}$ of the non-negative rational
functions. The article is based on the following methods:\\
1.We investigate non-negative trigonometrical polynomials and then
with the help of the Calderon transformation we proceed to the
power polynomials.\\
2.The way of constructing support hyperplanes to the convex sets
$Q_{P}$ and ${\sigma}_{P}$ is given in the paper.\\
Now we start with a more detailed description of the results of
this article.\\
Let us denote by ${\sigma}(N_{1},N_{2},N_{3})$ the set of the
trigonometrical polynomials
\begin{equation}
f(\alpha,\beta,\gamma)=\sum_{|k|{\leq}N_{1}}\sum_{|\ell|{\leq}N_{2}}\sum_{|m|{\leq}N_{3}}q(k,\ell,m)e^{i(k\alpha+\ell{\beta}+m{\gamma})}
\end{equation}
satisfying the condition
\begin {equation}
f(\alpha,\beta,\gamma){\geq}0,\quad
(\alpha=\bar{\alpha},\beta=\bar{\beta},\gamma=\bar{\gamma}).
\end{equation}
We denote by $Q(N_{1},N_{2},N_{3})$ the set of trigonometrical
polynomials of the class ${\sigma}(N_{1},N_{2},N_{3})$ admitting
representation
\begin{equation}
f(\alpha,\beta,\gamma)=\sum_{j=1}^{r}|F_{j}(\alpha,\beta,\gamma)|^{2}
\end{equation}
where
\begin{equation}
F_{j}(\alpha,\beta,\gamma)=\sum_{0{\leq}k{\leq}N_{1}}\sum_{0{\leq}\ell{\leq}N_{2}}\sum_{0{\leq}m{\leq}N_{3}}q_{j}(k,\ell,m)e^{i(k\alpha+\ell{\beta}+m{\gamma})}
\end{equation}
It is clear that the set $Q(N_{1},N_{2},N_{3})$ is convex.In this
article we give a method of constructing the support hyperplanes
to the set $Q(N_{1},N_{2},N_{3})$.Hence we receive a number of
general facts referring to the extremal points and faces of the
set $Q(N_{1},N_{2},N_{3})$.Here we also introduce concrete
examples of extremal points and faces.\\
Analogues results are received for the convex set
${\sigma}(N_{1},N_{2},N_{3})$ as well.In addition to the set of
the non-negative trigonometrical polynomials we shall introduce
the class ${\sigma}_{P}(2N_{1},2N_{2},2N_{3})$ of the power
non-negative  polynomials of the form
\begin{equation}
f(x,y,z)=\sum_{0{\leq}k{\leq}2N_{1}}\sum_{0{\leq}\ell{\leq}2N_{2}}\sum_{0{\leq}m{\leq}2N_{3}}a_{k,\ell,m}x^{k}y^{\ell}z^{m}
\end{equation}
By $Q_{P}(2N_{1},2N_{2},2N_{3})$ we denote the set of the power
non-negative  polynomials of the class
${\sigma}_{P}(2N_{1},2N_{2},2N_{3})$ admitting the representation
\begin{equation}
f(x,y,z)=\sum_{j=1}^{r}|F_{j}(x,y,z)|^{2}
\end{equation}
where $F_{j}(x,y,z) $are polynomials of x,y,z.With the help of the
Calderon transformation the results obtained for the classes of
the trigonometrical polynomials ${\sigma}(N_{1},N_{2},N_{3})$ and
$Q(N_{1},N_{2},N_{3})$ we transfer onto the classes of the power
polynomials ${\sigma}_{P}(2N_{1},2N_{2},2N_{3})$ and
$Q_{P}(2N_{1},2N_{2},2N_{3})$.Let us note that a number of
concrete examples of the extremal power polynomials is contained
in important works [2],[5].\\In our paper we consider the case of
three variables,but the obtained results can be easily transferred
to any number of variables.
\section{Main Notions}
Let S be a set of points of the space $R^{3}$.Let us denote by
${\Delta}=S-S$ the set of points $x{\in}R^{3}$ which can be
represented in the form $x=y-z,\quad y,z\in S$.The function
$\Phi(x)$ is called Hermitian positive on $\Delta$ if for any
points $x_{1},x_{2},...x_{N}{\in}S$ and numbers
${\xi}_{1},{\xi}_{2},...{\xi}_{N}$ the inequality
\begin{equation}
\sum_{i,j}{\xi}_{i}\bar{\xi}_{j}\Phi(x_{i}-x_{j}){\geq}0
\end{equation}
is true.We shall consider the lattice $S(N_{1},N_{2},N_{3})$
consisting of the points $M(k,\ell,m)$ where
$0\leq{k}{\leq}N_{1},0\leq{\ell}{\leq}N_{2},0{\leq}m{\leq}N_{3}$.
The set $\Delta(N_{1},N_{2},N_{3})$ consists of the points
$M(k,\ell,m)$ where $|k|{\leq }N_{1},|\ell|{\leq }N_{1},|m|{\leq}
N_{3}$.By $P(N_{1},N_{2},N_{3})$ we denote the class of functions
which are Hermitian positive on $\Delta(N_{1},N_{2},N_{3})$.\\
With each function $\Phi(k,\ell,m)$ from $P(N_{1},N_{2},N_{3})$\\
we associate the Toeplitz matrices (see [8],[9]:
\begin{equation}
  B(l,m) =
      \left[
\begin{array}{llll}
       \Phi (0,l,m)    & \Phi(1,l,m)      & \cdots & \Phi(N_1,l,m)    \\
       \Phi (-1,l,m)   & \Phi(0,l,m)      & \cdots & \Phi(N_1-1,l,m)  \\
       \hskip6mm\cdots & \hskip6mm\cdots  & \cdots & \hskip6mm\cdots  \\
       \Phi(-N,l,m)    & \Phi(-N_1+1,l,m) & \cdots & \Phi(0,l,m)
\end{array}
      \right]
   \label{3.3.1}
\end{equation}
>From the matrices $B(l,m)$ we construct the block Toeplitz
matrices
$$
  C_m =
       \left[
\begin{array}{llll}
      B(0,m)           & B(1,m)          & \cdots & B(N_2,m)            \\
      B(-1,m)          & B(0,m)          & \cdots & B(N_2-1,m)          \\
      \hskip6mm\cdots  &\hskip6mm\cdots  & \cdots & \hskip6mm\cdots     \\
      B(-N_2,m)        & B(-N_2+1,m)     & \cdots & B(0,m)
\end{array}
       \right] ,
\quad
 | m | \leq N_3
$$

       Finally from
$C_k$ we make yet another block Toeplitz matrix
\begin{equation}
    A(N_1,N_2,N_3) =
                     \left[
\begin{array}{llll}
       C_0               & C_1            & \cdots & C_{N_3}         \\
       C_{-1}            & C_0            & \cdots & C_{N_3-1}       \\
       \hskip2mm\cdots   &\hskip2mm\cdots & \cdots & \hskip2mm\cdots \\
       C_{-N_3}          & C_{-N_3+1}     & \cdots & C_0
\end{array}
             \right]
\label{313}
\end{equation}
\textbf{Proposition 1}(see [8],[9]).\emph{ Inequality (7) is
equivalent to the inequality\\ $A(N_1, N_2, N_3) \geq 0$}.\\
With the help of the function $\Phi(k,\ell,m)$\\
 we introduce the linear functional (see[7]):
\begin{equation}
L_{\Phi}(f)=\sum_{|k|{\leq}N_{1}}\sum_{|\ell|{\leq}N_{2}}\sum_{|m|{\leq}N_{3}}q(k,\ell,m)\Phi(k,\ell,m)
\end{equation}
\textbf{Proposition 2}(see [8],[9]).\emph{If\\
$\Phi(k,\ell,m){\in}P(N_1, N_2, N_3)$ and $f(\alpha,\beta,\gamma)=
|F(\alpha,\beta,\gamma)|^{2}$ then the relation
\begin{equation}
L_{\Phi}(f)=e^{\star}Ae{\geq}0
\end{equation}
holds.}\\
Here the matrix A is defined by relations (8),(9),the function
$F(\alpha,\beta,\gamma)$ and the vector e have the forms
\begin{equation}
F(\alpha,\beta,\gamma)=\sum_{0{\leq}k{\leq}N_{1}}\sum_{0{\leq}\ell{\leq}N_{2}}\sum_{0{\leq}m{\leq}N_{3}}q(k,\ell,m)e^{i(k\alpha+\ell{\beta}+m{\gamma})}
\end{equation}
where
\begin{equation}
e=col[h(0),h(1),...,h(N_{3}],
\end{equation}
\begin{equation}
h(m)=col[g(o,m),g(1,m),...,g(N_{2},m)],
\end{equation}
\begin{equation}
g(\ell,m)=col[d(0,\ell,m),d(1,\ell,m),...,d(N_{1},\ell,m)].
\end{equation}
\section{Support hyperplanes,extremal points and extremal faces of
 $Q(N_{1},N_{2},N_{3})$}
 Let $L$ be a linear a linear functional.The hyperplane
 $H=[L,\alpha]$ is said to bound the set $U$ if either \\
 $L(f){\geq}\alpha$ for all $f\in U$ or $L(f){\leq}\alpha$ for all $f\in U$.\\
 A hyperplane $H=[L,\alpha]$ is said to support a set U at a point
 $f_{0}{\in}U$ if $L(f_{0})=\alpha$ and if H bound U.\\
 Further we shall consider only such support hyperplanes $H$ which
 have at least one common point with the corresponding convex set
 U.\\
 A point $f_{0}$ in the convex set $U$ is called an extremal point
 of $U$ if there exists no non-degenerate line segment in $U$ that
 contains $f_{0}$ in its relative interior(see[10]).
 From Proposition 2 we obtain the following important assertion.\\
 \textbf{Corollary 1}.\emph{The set of support hyperplanes for
 $Q(N_{1},N_{2},N_{3})$ coincides with the set of hyperplanes
\begin{equation}
L_{\Phi}(f)=0
\end{equation}
where $\Phi(k,\ell,m){\in}P(N_{1},N_{2},N_{3})$ and the
corresponding matrix $A(N_{1},N_{2},N_{3})$ is such that
\begin{equation}
detA(N_{1},N_{2},N_{3})=0.
\end{equation}}
Let us denote by ${\nu}_{A}$ the dimension of the kernel of the
matrix $A(N_{1},N_{2},N_{3})$.If vector $e\ne0$ belongs to the
kernel of $A(N_{1},N_{2},N_{3})$ then the corresponding polynomial
$f(\alpha,\beta,\gamma)$ (see (1) and (14)-(18)) belongs to
$Q(N_{1},N_{2},N_{3})$ and satisfies relation (19).The convex hull
of such polynomials we denote by $D_{A}$.Using classical
properties of a convex set  we obtain the following
assertions. \\
\textbf{Corollary 2}.\emph{If ${\nu}_{A}=1$ and
$f(\alpha,\beta,\gamma){\in}D_{A}$ then the
$f(\alpha,\beta,\gamma)$ is the extremal polynomial in the class
$Q(N_{1},N_{2},N_{3})$.}\\
\textbf{Corollary 3}.\emph{If ${\nu}_{A}>1$ then $D_{A}$ is the
extremal face}. \\
\textbf{Corollary 4}.\emph{If $f(\alpha,\beta,\gamma)$ is an
extremal polynomial in the class $Q(N_{1},N_{2},N_{3})$ then there
exists a non trivial function
$\Phi(k,\ell,m){\in}P(N_{1},N_{2},N_{3})$ such that
$L_{\Phi}(f)=0$.}\\

 \VS
   {\bf Example 1 .}
   Let the relations
$$
  N_1 = N_2 = N_3 = 1
$$
be valid.  We set
$$
  \begin{array}{ll}
        B(0,0) =
                \left[
                     \begin{array}{cc}
                            1  &  0  \\
                            0  &  1
                     \end{array}
                \right] ,
&
        B(0,1) = B(1,0) =
                          \left[
                               \begin{array}{cc}
                                     0  &  0  \\
                                     0  &  0
                               \end{array}
                          \right]
\\
\,
\\
        B(1,1) =
                 \left[
                      \begin{array}{cc}
                            0  &  1  \\
                            1  &  0
                      \end{array}
                 \right] ,
&
        B(-1,1) =
                 \left[
                      \begin{array}{rr}
                            0  &  1  \\
                           -1  &  0
                      \end{array}
                 \right]
  \end{array}$$

      Then we have
\begin{equation}
      C_0 = E_4,
\quad
      C_1 =
           \left[
                \begin{array}{rrrr}
                       0  &  0  &  0  &  1  \\
                       0  &  0  &  1  &  0  \\
                       0  &  1  &  0  &  0  \\
                      -1  &  0  &  0  &  0
                \end{array}
           \right] ,
\quad
       A =
           \left[
                \begin{array}{ll}
                      E_4^{\,}  &  C_1^{\,}  \\
                      C_1^*     &  E_4
                \end{array}
           \right]
  \label{3.6.33}
\end{equation}
It is easy to see that the following linearly independent vectors

\begin{equation}
      \left\{
           \begin{array}{l}
                 e_1 = \mbox{col } [0,0,0,1,1,0,0,0]
\\
                 e_2 = \mbox{col } [0, 0, -1,0,0,1,0,0]
\\
                 e_3 = \mbox{col } [0,-1, 0,0,0,0,1,0]
\\
                 e_4 = \mbox{col } [-1,0,0,0,0,0,0,1]
           \end{array}
      \right.
  \label{3.6.35}
\end{equation}
form the basis of the kernel $A(1,1,1)$.The corresponding
polynomials $F_k (\alpha, \beta, \gamma)$ have the forms
\begin{equation}
      \left\{
            \begin{array}{ll}
                  F_1 (\alpha, \beta, \gamma) =
                  e^{i (\alpha + \beta)} + e^{i \gamma},
&
                  F_2 (\alpha, \beta, \gamma) =
                - e^{i \beta} + e^{i (\alpha + \gamma)}
\\
                  F_3 (\alpha, \beta, \gamma) =
                - e^{i \alpha} + e^{i (\beta + \gamma)},
&
                  F_4 (\alpha, \beta, \gamma) =
                - 1 + e^{i (\alpha + \beta + \gamma)}
            \end{array}
      \right.
  \label{3.6.36}
\end{equation}
It means that polynomials
\begin{equation}
f_k (\alpha, \beta, \gamma)=|F_k (\alpha, \beta,\gamma)|^{2},
\quad (1{\leq}k{\leq}4)
\end{equation}
belong to the face $D_{A}$.\\
 \textbf{Proposition 3}.\emph{The
polynomials $f_k (\alpha, \beta, \gamma)$ constructed by formulas
(20)and (21) are extremal in the class $Q(1,1,1)$.}\\
Proof.Let $B(0,0),B(0,1)$ and $B(1,0)$ be defined as above
 and let
\begin{equation}
 B(-1,1) =
                \left[
                     \begin{array}{cc}
                            0  &  {\gamma}_{1}  \\
                            -{\gamma}_{2}  &  0
                     \end{array}
                \right] ,
\end{equation}
\begin{equation}
         B(1,1) =
                          \left[
                               \begin{array}{cc}
                                     0  &  {\gamma}_{4}  \\
                                     {\gamma}_{3}  &  0
                               \end{array}
                          \right].
\end{equation}
We consider the cases when ${\gamma}_{j}=1,\quad |{\gamma}_{k}|<1,
\quad (k{\ne}j)$.In this cases vector $e_{j}$ (see(19))belongs to
the kernel of the corresponding matrix $A_{j}(1,1,1)$ and\\
${\nu}_{j}=1$.Hence the polynomials $f_{j}(\alpha,\beta,\gamma)$
are extremal.\\
\section{Support hyperplanes,extremal points and extremal faces of
 $\sigma(N_{1},N_{2},N_{3})$}
We will say that the function
$\Phi(k,\ell,m){\in}P(N_{1},N_{2},N_{3})$ is extendible if
$\Phi(k,\ell,m)$ admits an extension to the function of the class
$P(\infty,\infty,\infty)$.
We will use the following Rudin's result [7].\\
\textbf{Proposition 4}.\emph{A function
 can be extended to a
member of $P(\infty,\infty,\infty)$ if and only if
$L_{\Phi}(f){\geq}0$ for every
$f{\in}\sigma(N_{1},N_{2},N_{3})$}.\\
>From Proposition 4 we deduce the following important assertion.\\
\textbf{Corollary 5.}\emph{The set of support hyperplanes of $
\sigma(N_{1},N_{2},N_{3})$ coincides with the set of hyperplanes
$L_{\Phi}(f)=0$ where $\Phi(k,\ell,m)$ is an extendible function
from $P(N_{1},N_{2},N_{3})$ and there exists a non-trivial
polynomial $f_{0}{\in}\sigma(N_{1},N_{2},N_{3})$ such that}\\
$L_{\Phi}(f_0)=0$.\\
The convex hull of polynomials
$f_{0}{\in}\sigma(N_{1},N_{2},N_{3})$ satisfying relation\\
$L_{\Phi}(f_0)=0$ we denote by $D_{\Phi}$.The number of linearly
independent polynomials $f_0$ from $D_{\Phi}$ we denote by
${\nu}_{\Phi}$.From classical properties of a convex set  we
obtain the following assertions.\\
 \textbf{Corollary 6}.\emph{1.If ${\nu}_{\Phi}=1$ and $f_0{\in}D_{\Phi}$
then $f_0$ is the extremal polynomial in the
class  $\sigma(N_{1},N_{2},N_{3})$.\\
2.If ${\nu}_{\Phi}>1$ then $D_{\Phi}$ is the
extremal face in the class $\sigma(N_{1},N_{2},N_{3})$.\\
3.If $f_0$ is an extremal polynomial in the class
$\sigma(N_{1},N_{2},N_{3})$ then there exists
 a non-trivial
extendible function $\Phi(k,\ell,m){\in}P(N_{1},N_{2},N_{3})$ such
that \\$L_{\Phi}(f_0)=0$}.\\
If $\Phi(k,\ell,m)$ is an extendible function then there exists a
positive measure $\mu(\alpha,\beta,\gamma)$ such that (Bochner
theorem)
 \begin{equation}
\Phi(k,\ell,m)=\frac{1}{(2\pi)^{3}}\int_{G}e^{i(k\alpha+\ell{\beta}+m\gamma)}d\mu
\end{equation}
where Domain $G$ is defined by the inequalities
$-\pi{\leq}\alpha,\beta,\gamma{\leq}\pi$.From relations (1),(10)
and (24)we deduce the following well-known representation (see[7])
\begin{equation}
L_{\Phi}(f)=\frac{1}{(2\pi)^{3}}\int_{G}f(\alpha,\beta,\gamma)d\mu.
\end{equation}
\textbf{Corollary 7}.\emph{If $f_{0}(\alpha,\beta,\gamma)$ belongs
to the class $\sigma(N_{1},N_{2},N_{3})$ and in a certain point
$f_{0}({\alpha}_0,{\beta}_0,{\gamma}_0)=0$ then
$f_{0}(\alpha,\beta,\gamma)$ is either extremal or belongs to the
extremal face.}\\
\section{The power non-negative polynomials}
We shall use the linear Calderon transformation
\begin{equation}
(Cf)(\alpha,\beta,\gamma)=f_{1}(x,y,z)
\end{equation}
which is defined by formulas \\
$e^{i\alpha}=\frac{x+i}{x-i},\quad
e^{i\beta}=\frac{y+i}{y-i},\quad \frac{y+i}{y-i}$,
\begin{equation}
f_{1}(x,y,z)=f(\frac{x+i}{x-i},\frac{y+i}{y-i},\frac{y+i}{y-i})(x^{2}+1)^{N_1}(y^{2}+1)^{N_2}(z^{2}+1)^{N_3}.
\end{equation}
\textbf{Proposition 5}(see [9],Ch.3).\emph{The Calderon
transformation C maps $\sigma(N_{1},N_{2},N_{3})$ onto
${\sigma}_{P}(2N_{1},2N_{2},2N_{3})$ and $Q(N_{1},N_{2},N_{3})$
onto$Q_{P}(2N_{1},2N_{2},2N_{3})$}.\\
Using linearity of the operator C we deduce from Proposition 5 the
following assertion.\\
\textbf{Corollary 8}.\emph{The Calderon transformation C maps the
extremal points and faces of $\sigma(N_{1},N_{2},N_{3})$ and
$Q(N_{1},N_{2},N_{3})$ onto extremal points and faces of
$\sigma_{P}(2N_{1},2N_{2},2N_{3})$ and
$Q_{P}(2N_{1},2N_{2},2N_{3})$ respectively.}\\
\textbf{Example 2}.Let us consider the polynomials \\
\begin{equation}
P_{1}(x,y,z)=(xyz-z+y+x)^{2},
\end{equation}
\begin{equation}
P_{2}(x,y,z)=(yz-xz+xy+1)^{2},
\end{equation}
\begin{equation}
P_{3}(x,y,z)=(yz-xz-xy-1)^{2},
\end{equation}
\begin{equation}
P_{4}(x,y,z)=(xz+yz+xy-1)^{2},
\end{equation}
\textbf{Proposition 6}.\emph{Polynomials $P_{k}(x,y,z)$
(k=1,2,3,4) are extremal polynomials in the classes $Q_{P}(2,2,2)$
and
$\sigma_{P}(2,2,2)$.}\\
Proof.Using Example 1 and Corollary 8 we deduce that polynomials
$P_{k}(x,y,z)$ (k=1,2,3,4) are extremal in the class
$Q_{P}(2,2,2)$. We remark that\\
 $deg P_{k}(x,y,z){\leq}6$. It is
proved for this case (see[1])that the extremal polynomials in the
class $Q_{P}$ are extremal in the class $\sigma_{P}(2,2,2)$
as well.The proposition is proved.\\
We remark that an extremal
polynomial in the class $Q_{P}$ is given in the paper [1].
\section{Extremal trigonometrical and power polynomials in the
classes $\sigma$ and $\sigma_{P}$} In my book [9] I give the
method of constructing trigonometrical polynomials belonging to
$\sigma(N_{1},N_{2},N_{3})$ but not belonging to
$Q(N_{1},N_{2},N_{3})$.This method can also be used for
constructing extremal polynomials in the classes
$\sigma(N_{1},N_{2},N_{3})$ and $\sigma_{P}(N_{1},N_{2},N_{3})$.
For illustrating this fact we consider the following example.\\
\textbf{Example 3}.Let us introduce the polynomial
\begin{equation}
f_{0}({\alpha},{\beta},{\gamma})=4-cos({\alpha}+{\beta}+{\gamma})-
cos(-{\alpha}+{\beta}+{\gamma})-cos({\alpha}-{\beta}+{\gamma})+cos({\alpha}+{\beta}-{\gamma})
\end{equation}
It is shown in the book ([9]Ch.3)that
$f_{0}({\alpha},{\beta},{\gamma})-m (0<m{\leq}4-2^{3/2})$ belongs
to
$\sigma(1,1,1)$ but does not belongs to $Q(1,1,1)$.\\
Now we consider the polynomial
\begin{equation}
f({\alpha},{\beta},{\gamma})=2^{3/2}-cos({\alpha}+{\beta}+{\gamma})-
cos(-{\alpha}+{\beta}+{\gamma})-cos({\alpha}-{\beta}+{\gamma})+cos({\alpha}+{\beta}-{\gamma})
\end{equation}
The last equality we rewrite in the form
\begin{equation}
f({\alpha},{\beta},{\gamma})=2[2^{1/2}-{cos\alpha}cos({\beta}+{\gamma})+
{sin\alpha}sin({\beta}-{\gamma})]
\end{equation}
\textbf{Proposition 7}.\emph{The polynomial
$f({\alpha},{\beta},{\gamma})$ is an extremal one in the class
$\sigma(1,1,1)$.}\\
Proof.It follows from formula (34) that the polynomial
$f({\alpha},{\beta},{\gamma})$ has the following zeroes:\\
1. $\alpha_{1}={\pi}/4,\beta(0,1),\gamma(0,1)$.\\
2. $\alpha_{2}=-{\pi}/4,\beta(0,0),\gamma(0,0)$.\\
3. $\alpha_{3}={\pi}/4,\beta(0,-1),\gamma(0,-1)$\\
4. $\alpha_{4}=-{\pi}/4,\beta(0,-2),\gamma(0,-2)$.\\
5. $\alpha_{5}=3{\pi}/4,\beta(-1,-1),\gamma(-1,-1)$.\\
6. $\alpha_{6}=-3{\pi}/4,\beta(-1,0),\gamma(-1,0)$.\\
7. $\alpha_{7}=3{\pi}/4,\beta(1,-1),\gamma(1,-1)$.\\
8. $\alpha_{8}=-3{\pi}/4,\beta(1,0),\gamma(1,0)$.\\
Since the polynomial $f({\alpha},{\beta},{\gamma})$ is
non-negative all its first derivatives in the points of zeroes are
equal to zero.This is also true for the polynomials
$g({\alpha},{\beta},{\gamma}){\in}\sigma(1,1,1)$ and such that
$g({\alpha},{\beta},{\gamma}){\leq}f({\alpha},{\beta},{\gamma})$.
Thus we obtain 32 linear equations on 27 coefficients of the
polynomial $g({\alpha},{\beta},{\gamma})$.This linear system is of
26th rank,that is only $cf (c=const)$ satisfies this system.The
proposition is proved.\\
>From Proposition 7 using the Calderon transformation we obtain the
following assertion.\\
\textbf{Proposition 8.}\emph{The polynomial
\begin{equation}
f(x,y,z)=2^{3/2}(1+x^{2})(1+y^{2})(1+z^{2})+8z(y+x)(yx-1)-2(z^{2}-1)[(yx+1)^{2}-(x-y)^{2}]
\end{equation}
is an extremal in the class $\sigma_{P}(2,2,2$)}.\\
As it was mentioned in the introductory part  some other extremal
polynomials in the class $\sigma_{P}$were known earlier [2].
\section{Non-negative rational functions}
Let us consider the class ${\sigma}_{R}(N,M)$ of the non-negative
rational functions of the form \\
$R(x,y)=\frac {p(x,y)}{q(x,y)}$\\
 where $p(x,y)$ and $q(x,y)$
are polynomials with real coefficients and $deg p(x,y,z){\leq}N,$
and $deg q(x,y,z){\leq}M$.\\
 As Artin's result
 shows the function $R(x,y,z)$ can be represented in the form
\begin{equation}
R(x,y)=\sum_{k=1}^{r}\frac {p_{k}^{2}(x,y)}{q_{k}^{2}(x,y)}
\end{equation}
where $p_{k}(x,y,z)$ and $q_{k}(x,y,z)$ are polynomials.It is
clear that  ${\sigma}_{R}(N,M)$ is a convex set.\\
\textbf{Proposition 9}.\emph{If $R(x,y,z)$ is an extremal function
of the convex set $\sigma(N,M)$ then $R(x,y,z)$ admits the
representation \\
\begin{equation}
R(x,y)=\frac {p^{2}(x,y)}{q^{2}(x,y)}
\end{equation}
 where $p(x,y,z)$ and
$q(x,y,z)$ are polynomials and $p(x,y)^{2}$
is extremal in the class of the polynomials $Q_{P}$.}\\
Thus the results concerning the class $Q_{P}$(section 4) can be
useful for the investigation of the non-negative rational
functions.\\
I am very grateful to professor B.Reznick for having sent me a
number of his papers which stimulated my work.I greatly thank
doctor I.Tydniouk for his help in counting the rank of the linear
system in Example 3.

\begin{center}\textbf{References}\end{center}
1.\textbf{M.D.Choi,M.Knebush,T.Y.Lam and
B.Reznick},\emph{Transversal zeros and positive semidefinite
forms},Lecture Notes in Math.959,Springer-Verlag (1982)273-298.\\
2.\textbf{M.D.Choi and T.Y.Lam},\emph{Extremal positive
semidefinite forms},Math.Ann.231 (1977)1-18\\
3.\textbf{D.Hilbert},\emph{Uber die Darstellung definiter Formen
als Summe von Formenquadraten},Math.Ann.32 (1888)342-350.\\
4.\textbf{T.S.Motzkin},\emph{Selected papers},Birkhauser,1983.\\
5.\textbf{B.Reznick},\emph{Some Concrete Aspects of Hilbert's 17th
Problem},Contemporary Mathematics,253,(2000),251-272.\\
6.\textbf{R.M.Robinson},\emph{Some definite polynomials which are
not sums of squares of real polynomials},Selected questions of
algebra and logic,Nauka,Novosibirsk (1973)264-282.\\
7.\textbf{W.Rudin},\emph{The Extension Problem for Positive
Definite Functions},Ilinois J.Mat.(1963)532-539.\\
8.\textbf{L.A.Sakhnovich},\emph{Effective construction of
non-extendible Hermitian-positive functions of several
variables},Funct.Anal.Appl.,14,No 4(1980),55-60.\\
9.\textbf{L.A.Sakhnovich},\emph{Interpolation Theory and its
Applications},Kluwer Acad.Publishers,(1997).\\
10.\textbf{S.R.Lay}\emph{Convex Sets and their
Applications},A.Wiley Interscience Publication,(1982).

\end{document}